\documentclass[12pt,reqno]{amsart}
\usepackage{amsmath, caption, floatflt, graphicx, verbatim, array, longtable, hhline,
moreverb}
\newtheorem{theorem}{Theorem}

\theoremstyle{definition}

\theoremstyle{remark}

\newcommand{\Area}{\operatorname{Area}}
\newcommand{\length}{\operatorname{length}}
\newcommand{\kk}{\mathbf k}

\theoremstyle{plain}

\newcommand{\R}{{\mathbf{R}}}
\newcommand{\C}{{\mathcal{C}}}
\newcommand{\N}{{\mathbf{N}}}
\newcommand{\T}{\mathcal{T}}

\newcommand{\D}{\mathcal{D}}

\newcommand{\I}{\mathcal{I}}
\newcommand{\J}{\mathcal{J}}

\def\NN{\mathbf{N}}

\newcommand{\OO}{{\mathfrak{O}}}

\newcommand{\F}{{\mathfrak{F}}}
\newcommand{\FF}{{\mathfrak{F}}}
\newcommand{\FQ}{{\mathfrak{F}_{_Q}}}

\newcommand{\ga}{\gamma}

\newcolumntype{C}{>{$}c<{$}}
\newcolumntype{D}{>{$} c<{$}}
\newcolumntype{L}{>{$}l<{$}}
\newcolumntype{R}{>{$}r<{$}}
\def\s{_{_{_{}}}} 

\begin{document}

\title[On the Intervals of a Third between Farey Fractions]
{On the Intervals of a Third between Farey Fractions}

\author[Cobeli, Zaharescu]{Cristian Cobeli, Alexandru Zaharescu}\footnotetext{CC is partially supported by the CERES Programme of the Romanian Ministry
of Education and Research, contract 4-147/2004.}

\address{CC: Institute of Mathematics of the
Romanian Academy, P.O. Box 1-764,
Bucharest 70700, Romania}
\email{cristian.cobeli@imar.ro}

\address{AZ: Institute of Mathematics of the
Romanian Academy, P.O. Box 1-764,
Bucharest 70700, Romania}

\address{University of Illinois at Urbana Champaign, Urbana, IL}
\email{zaharesc@math.uiuc.edu}

\thanks{}

\begin{abstract}

The spacing distribution between Farey points has drawn attention in recent years.
It was found that the gaps $\ga_{j+1}-\ga_j$ between consecutive 
elements of the Farey sequence produce, as $Q\to\infty$, a limiting measure.
Numerical computations suggest that for any $d\ge 2$, the gaps 
$\ga_{j+d}-\ga_j$ also produce a limiting measure whose support is
distinguished by remarkable topological features.
Here we prove the existence of the spacing distribution for $d=2$ and
characterize completely the corresponding support of the measure.
\end{abstract}
\maketitle

\section{Introduction}

Let $\FQ=\{\ga_1,\dots ,\ga_N\}$ be the Farey sequence of order $Q$,
which is defined to be the set of all subunitary irreducible fractions
with denominators $\le Q$, arranged in ascending order.
For any interval $\I\subset [0,1]$, we write $\F_Q(\I)=\F_Q\cap\I$.
The cardinality of $\F_Q(\I)$ is well known to be 
$N_\I(Q)=3|\I|Q^2/\pi^2+O(Q\log Q).$ When $\I=[0,1]$ we write shortly
$N(Q)$ instead of $N_{[0,1]}(Q)$.
Since $\FQ$ contains a large number of fractions obtained by a combined
process of division, sieving and sorting of integers from $[1,Q]$, one
would apriori expect little or even no special structure in the set of all
differences between consecutive fractions 
(which we also call {\it intervals of a second}). 
Though, this expectation is not
fulfilled. This is sustained from many points of view by a
series of authors, such as 
Franel \cite {Franel}, 
Kanemitsu, Sita Rama Chandra Rao and Siva Rama Sarma \cite {KSS}, 
Hall and Tenenbaum \cite {HT1}, \cite {HT2}, 
Hall \cite {Hall}, 
Augustin, Boca and the authors~\cite{ABCZ}, 
who have studied the set of gaps between consecutive Farey fractions.
A  regularity is expected also in the set of larger gaps
$\ga^{(d+1)}-\ga'$, where $\ga'$ runs over $\{\ga_1,\dots ,\ga_{N-d}\}$ and
$d\ge 2$. (We use up-scripts, such as $\ga',\ga'',\ga''',\dots $ to write consecutive
elements of $\FQ$.) It is our object to treat here the case $d=2$, that
is, the case of {\it intervals of a third}.

\begin{center}
\begin{figure}[h]
    \begin{minipage}[b]{\linewidth}
      \centering
      \includegraphics*[width=0.8\linewidth,   bb=  -18 -29 700 700
       ]{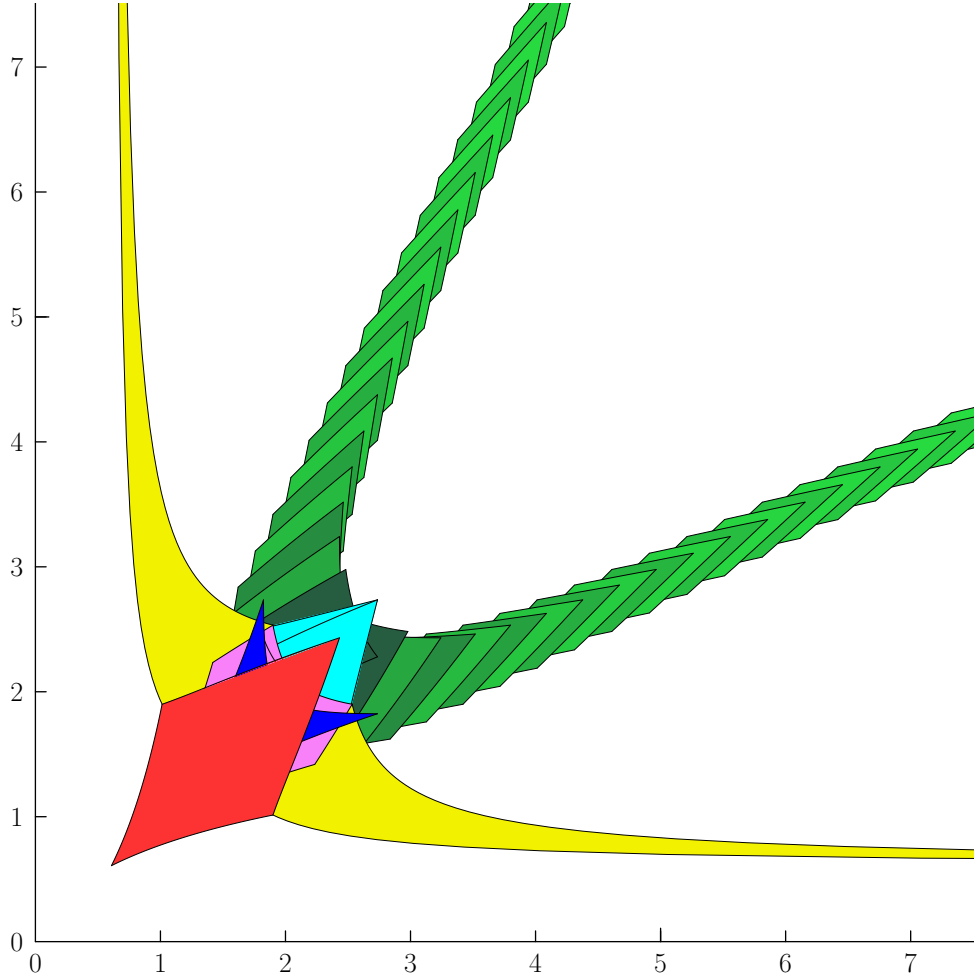}
          \caption{
The support of $\mu_{2,2}$.
                   }\label{Figure1}
     \end{minipage}
\end{figure}
\end{center}

Geometrical representations of the set of pairs of neighbor intervals
of fractions from $\FQ$ created for different values of $Q$ 
reveals sets of points whose density concentrates on different parts
of the plane. The aesthetical qualities of the pictures catches
attention immediately. For for any $d\ge 1$ they look like a  
swallow and the main topological distinctions are in the number of folds of the
tail. Thus, when $d=1$ (neighbor pairs of intervals of a second) the
swallow has a one-fold tail (see \cite{ABCZ}).
When $d=2$, the case treated in the present paper, the swallow has a two-fold tail 
(see Figure~\ref{Figure1}) and in Section~\ref{secSupport} we have calculated
explicitely the equations of the frontier.
In the cases  $d\ge 3$ the tail appears always to have a
three-folded tail, but this is more complex and its characterization
will appear in a separate paper.

Given $N$ real numbers $x_1 \le x_2 \le\dots \le x_N$ with mean
spacing $1$, we consider \emph{the $h$-th level of intervals of a third
probability} $\mu_{2,h}$  on $ R_+^h $, defined, for $f\in C_c([0,\infty))$, by
	\begin{equation*}
		\int_{[0,\infty)^h} f d\mu_{2,h} = 
		\frac 1{N-h-1} 
		\sum_{j=1}^{N-h-1} f(x_{j+2}-x_j, x_{j+3} - x_{j+1}, \dots ,
				x_{j+h+1}-x_{j+h-1})\,.
	\end{equation*}

In our case, we normalize $\F_Q(\I)$ to get
the sequence $x_j=N(Q,\I)\ga_j/|\I|$, $1\le j\le N(Q,\I)$ with mean spacing
equal to one. Accordingly, we get the sequence $(\mu_{2,h}^{Q,\I})_{Q\ge 1}$
of the $h$-th level of intervals of a third probabilities on $[0,\infty)^h$.
We show that this sequence converges, as $Q\to\infty$, to a probability 
measure $\mu_{2,h}$, which is independent of $\I$,
and  can be expressed explicitly.

For any  $\gamma_i=a_i/q_i$ and $\gamma_j=a_j/q_j$ in $\FQ$, we set
$\Delta(\gamma_i,\gamma_j)=\Delta(i,j)=-
\bigl|
\begin{smallmatrix}
                       a_i&a_j \\ q_i&q_j
\end{smallmatrix}
\bigr|\,.
$
This is the numerator of the difference $\ga_j-\ga_i$.
It is well known that $\Delta(\ga',\ga'')=-1$ for any consecutive elements
of $\FQ$, and it turns out that this equality is responsible for the
existence of the $h$-spacing distribution of the Farey sequence.
Though, this relation is no longer true for larger intervals,
but there is a convenient replacement. To see this, let us note that
a Farey fraction can be uniquely determined by its two
predecessors. Indeed, if $\frac {a'}{q'}<\frac{a''}{q''}<\frac{a'''}{q'''}$
are consecutive fractions of $\FQ$, we have
$ a'''=ka''-a'$ and $q'''=kq''-q'$,
where $k=\Delta(\ga',\ga''')=\big[\frac{q'+Q}{q''}\big]$.  

The basic idea of our procedure is to parametrize the set of $h$-tuples of
intervals of a third in terms of just two variables that run over a
completely described domain.
The set of pairs of consecutive denominators of fractions in 
$\FQ$ are exactly the elements of
  \begin{equation*}
    \begin{split}
     \big\{(q',q'')\colon\ 1\le q',q''\le Q,\ q'+q''>Q \ 
           \text{and}\ (q',q'')=1\big\}\,.
    \end{split}
  \end{equation*}
Since we are mainly interested in what happens when $Q\rightarrow\infty$, 
we reduce the scale $Q$ times, and consider the
background triangle $\T=\big\{(x,y)\colon\ 0<x\le 1,\ x+y>1 \big\}$,
called \emph{the Farey triangle}. We split it into a series of polygons as follows.
Firstly, for each $(x,y)\in \R^2$, we set $L_0(x,y)=x$, $L_1(x,y)=y$, and then,
for $i\geq 2$, we define recursively:
	\begin{equation*}\label{eq2}
		L_i(x,y)=\left[ \frac{1+L_{i-2}(x,y)}{L_{i-1}(x,y)} \right]
		L_{i-1}(x,y)-L_{i-2}(x,y)\,. 
	\end{equation*}
Then, as in \cite{BCZ1}, we consider the map
	\begin{equation*}
		{\bf k} :\T \rightarrow (\N^*)^{h} ,\qquad
		{\bf k}(x,y)=\big( k_1(x,y),\dots ,k_{h} (x,y)\big),
	\end{equation*}
where
$k_i(x,y)=\left[ \frac{1+L_{i-1}(x,y)}{L_i(x,y)} \right]$.
The functions $k_i(x,y)$ are locally constant, and the subsets of $\T$ on
which they are constant plays a special role. Thus, for any
${\bf k} \in (\N^*)^{h}$, we get the convex polygon
	\begin{equation*}
		\T_{\bf k} =\big\{ (x,y)\in \T \colon\
		{\bf k}(x,y)={\bf k} \big\} \,.
	\end{equation*}
Notice that
$\T =\bigcup\limits_{{\bf k}\in (\N^*)^{h}} \T_{\bf k}$
and
$\T_{\bf k} \cap \T_{{\bf k}^\prime} =\emptyset$
whenever ${\bf k} \neq {\bf k}^\prime$.

Next we consider the application $\Phi_{2,h}\colon\ \T\to (0,\infty)^h$
defined by
$$\Phi_{2,h}(x,y)=\frac 3{\pi^2}\left(
    \frac{k_1(x,y)}{L_0(x,y)L_2(x,y)},
    \frac{k_2(x,y)}{L_1(x,y)L_3(x,y)},\dots,
    \frac{k_h(x,y)}{L_{h-1}(x,y)L_{h+1}(x,y)}
    \right).
$$
Our main result shows that, indeed, for $Q\to\infty$, the sequence
$(\mu_{2,h}^{Q,\I})_{Q\ge 1}$ converges to a measure and
$\Phi_{2,h}(x,y)$ is the needed tool to describe its support.
\begin{theorem}\label{Theorem1}
The sequence $(\mu_{2,h}^{Q,\I})_{Q\ge 1}$ converges weakly to a probability
measure $\mu_{2,h}$, which is independent of $\I$. 
The support $\D_{2,h}$ of $\mu_{2,h}$ is the closure of the range of 
$\Phi_{2,h}$, and 
  \begin{equation*}
    \mu_{2,h}(\C)=2\mathrm{Area}(\Phi_{2,h}^{-1}(\C))\,,
  \end{equation*}
for any parallelepiped $\C=\prod_{j=1}^h (\alpha_j,\beta_j)\subset (0,\infty)^h$.  
\end{theorem}

In Table~\ref{Table1} from Section~\ref{secSupport} we provide explicit formulae 
for all the pieces that form $\D_{2,2}$.

\section{The Existence of the Limiting Measure}\label{secMeasure}

It is plain that in order to prove Theorem~\ref{Theorem1}, it suffices
to see the effect of $\mu_{2,h}$ on bounded parallelepipeds. 
For any $\C=\prod_{j=1}^h (\alpha_j,\beta_j)\subset (0,\infty)^h$, we define
  \begin{equation*}
    \mu_{2,h}^{Q,\I}(\C):=
	\frac{1}{N_\I(Q)} \ \cdot\# 
	\left\{ \gamma_j \in \FF_\I(Q)\colon\
	\frac{\alpha_i \vert \I\vert}{N_\I(Q)} < 
	\gamma_{j+i+1}-\gamma_{j+i-1}
	<\frac{\beta_i \vert \I\vert}{N_\I(Q)},\ i=1,\dots,h \right\} .
  \end{equation*}
We have to show that the sequence $\{\mu_{2,h}^{Q,\I}\}_Q$ is
convergent when $Q\to\infty$ and the limit is independent of $\I$.
In the beginning we treat the case of the complete interval $\I=[0,1]$.
\subsection{The case $\I=[0,1]$}
In the following we write shortly
$\mu_{2,h}^{Q}$ instead of $\mu_{2,h}^{Q,[0,1]}$.

With the notations from the Introduction, we see that
$\gamma_{j+i+1}-\gamma_{j+i-1}=k_{j+i}/q_{j+i+1}q_{j+i-1}$. Then
$\mu_{2,h}^{Q}(\C)$ can be written as
  \begin{equation}\label{eqdoi}
    \mu_{2,h}^{Q}(\C)=
	\frac{1}{N(Q)} \ \cdot \# 
	\left\{ \gamma_j \in \FF(Q)\colon\
	\frac{N(Q)}{\beta_i} < 
	\frac{q_{j+i+1}q_{j+i-1}}{k_{j+i}}
	<\frac{N(Q)}{\alpha_i},\ i=1,\dots,h \right\}. 
  \end{equation}
Knowing that $q_{j+i}=QL_i(q_j/Q,q_{j+1}/Q)$, we consider the set
  \begin{equation}\label{eqOQ}
    \Omega^{Q}(\C)=
	\Big\{ (x,y)\in Q\T\colon\
	\frac{N(Q)}{Q^2\beta_i} < 
		\frac{L_{i-1}\big(\frac xQ,\frac yQ\big)
			L_{i+1}\big(\frac xQ,\frac yQ\big)}
		     {k_{i}\big(\frac xQ,\frac yQ\big)}
	<\frac{N(Q)}{Q^2\alpha_i},\ i=1,\dots,h \Big\}. 
  \end{equation}
Since neighbor denominators in $\FF_Q$ are always coprime, relation
\eqref{eqdoi} turns into 
  \begin{equation*}
    \mu_{2,h}^{Q}(\C)=
	\frac{1}{N(Q)} \ \cdot \# 
	\left\{ (x,y)\in \Omega^{Q}(\C)\cap \NN^2\colon\
		\gcd(x,y)=1\right\}. 
  \end{equation*}

Next, we select the points with coprime coordinates using M\"obius
summation (cf.~\cite[Lemma~2]{ABCZ}), and we find that  
  \begin{equation}\label{eqtrei}
    \mu_{2,h}^{Q}(\C)=
	\frac{1}{N(Q)} 
	\Big( \frac{6}{\pi^2}\Area\big(\Omega^{Q}(\C)\big)
		+O\big(\length\big(\partial\Omega^{Q}(\C)\big)\log Q\big)\Big). 
  \end{equation}
Splitting $\T$ into the series of polygons $\T_\kk$, we see that the
error term in \eqref{eqtrei} is $O(Q\log Q)$. In the main term, we
replace $\Omega^{Q}(\C)$ by the bounded set
$\Omega(\C)=\Omega^{Q}(\C)/Q$. These yield
  \begin{equation}\label{eqpatru}
    \mu_{2,h}^{Q}(\C)=
	\frac{6Q^2}{\pi^2N(Q)} \Area\big(\Omega(\C)\big)
		+O_\C\Big(\frac{\log Q}{Q}\Big). 
  \end{equation}
It remains to replace in \eqref{eqpatru} the set $\Omega(\C)$ by a set
as in \eqref{eqOQ}, but with bounds independent of $Q$ in the corresponding
inequalities. This set is
  \begin{equation}\label{eqOOQ}
    \OO(\C):=
	\left\{ (x,y)\in Q\T\colon\
	\frac{3}{\pi^2\beta_i} < 
		\frac{L_{i-1}(x,y)L_{i+1}(x,y)}{k_{i}(x,y)}
	<\frac{3}{\pi^2\alpha_i},\ i=1,\dots,h \right\}. 
  \end{equation}
Notice that $\OO(\C)$ is exactly $\Phi^{-1}_{2,h}(\C)$.
The replacement does not change the error term because, via 
$N(Q)=3Q^2/\pi^2+O(Q\log Q)$, we have:
    \begin{equation}\label{eqmax}
	 \max\limits_{1\leq i\leq h} 
	\Big\{ \Big\vert \frac{N(Q)}{\alpha_iQ^2}-\frac{3}{\pi^2\alpha_i} \Big\vert, 
	\Big\vert \frac{N(Q)}{\beta_iQ^2}-\frac{3}{\pi^2\beta_i} \Big\vert
		\Big\}= O_\C\Big( \frac{\log Q}{Q}\Big) \, ,
   \end{equation}
which implies
    \begin{equation}\label{eqer1}
		\Area \big( \Omega({\C}) \triangle \OO(\C)\big)
		= O_\C\Big( \frac{\log Q}{Q}\Big) \, .
   \end{equation}
Therefore, by \eqref{eqOOQ} and \eqref{eqer1}, we get
  \begin{equation}\label{eqfizu}
    \mu_{2,h}^{Q}(\C)=2 \Area\big(\OO(\C)\big)
		+O_\C\Big(\frac{\log Q}{Q}\Big). 
  \end{equation}

In particular, this gives 
$\mu_{2,h}(\C)=\lim_{Q\to\infty}\mu_{2,h}^{Q}(\C)=2
\Area\big(\OO(\C)\big)$, 
concluding the proof of the theorem when $\I=[0,1]$.

\subsection{The short interval case}
Suppose now that $\I\subset [0,1]$ is fixed. In order to impose the
condition that only the fractions from $\I$ are involved in the
calculations, we employ the fundamental property of neighbor fractions
in $\FF_Q$. This says that if $\ga'=a'/q'$ and $\ga''=a''/q''$ are
consecutive then $a''q'-a'q''=1$. Consequently, $a''\equiv (q')^{-1}$
(mod $q''$), and this allows us to write the fraction $a''/q''$ in
terms of $q'$ and $q''$. Thus 
	\begin{equation*}
	a''/q''\in\I \iff
		(q')^{-1} \pmod {q''}\in q''\I\,.
	\end{equation*}  

This time we have to estimate
  \begin{equation}\label{equnspe}
    \mu_{2,h}^{Q,\I}(\C)=
	\frac{1}{N_\I(Q)} \ \cdot \# 
	\left\{ (q',q'')\in Q\T\colon\
		\begin{array}{l}
	\frac{N_\I(Q)}{|\I|Q^2\beta_i} < 
		\frac{L_{i-1}\big(\frac {q'}Q,\frac {q''}Q\big)
			L_{i+1}\big(\frac {q'}Q,\frac {q''}Q\big)}
		     {k_{i}\big(\frac {q'}Q,\frac {q''}Q\big)}
	<\frac{N_\I(Q)}{|\I|Q^2\alpha_i},\\ 
			\text{for}\ i=1,\dots,h;\quad 	(q')^{-1} \pmod {q''}\in q''\I
		\end{array}
	\right\}. 
  \end{equation}
We may write \eqref{equnspe} as
  \begin{equation}\label{eqdoispe}
    \mu_{2,h}^{Q,\I}(\C)=
	\frac{1}{N_\I(Q)} 
		\sum_{q=1}^Q N_q(\J^Q_\C(q),q\I)\,,
  \end{equation}
where
  \begin{equation*}
	N_q(\J_1,\J_2)=\#\big\{
		(m,n)\in\J_1\times\J_2\colon\ mn\equiv 1 \pmod q
	\big\}\,,
  \end{equation*}
for any $\J_1,\J_2\subset [0,Q-1]$ and
  \begin{equation*}
	\J^Q_\C(q)=
	\left\{ x\in(Q-q,Q]\colon\
		\begin{array}{l}
	\frac{N_\I(Q)}{|\I|Q^2\beta_i} < 
		\frac{L_{i-1}\big(\frac {q'}Q,\frac {q''}Q\big)
			L_{i+1}\big(\frac {q'}Q,\frac {q''}Q\big)}
		     {k_{i}\big(\frac {q'}Q,\frac {q''}Q\big)}
	<\frac{N_\I(Q)}{|\I|Q^2\alpha_i},\ i=1,\dots,h
		\end{array}
	\right\}. 
  \end{equation*}
For the best available technique to estimate the size of
$N_q(\J_1,\J_2)$ one requires bounds for Kloosterman sums
(cf.~\cite{BCZ}).  This is done when $\J_1$ and $\J_2$ are intervals,
but it may be easily extended for finite unions of subintervals of
$[0,q-1]$ (as the set $\J^Q_\C(q)$ is), even with the same same
formula. For our needs here, it suffices a version with a slightly
weaker term: 
  \begin{equation}\label{eqtreispe}
	N_q(\J^Q_\C(q),q\I)=\frac{\varphi(q)|\J^Q_\C(q)|\cdot|\I|}{q}
	+O_{\C,\varepsilon}\big(q^{1/2+\varepsilon}\big)\,.
  \end{equation}
Inserting \eqref{eqtreispe} into \eqref{eqdoispe}, we get
  \begin{equation}\label{eqpaispe}
    \mu_{2,h}^{Q,\I}(\C)=
	\frac{|\I|}{N_\I(Q)} 
		\sum_{q=1}^Q \frac{\varphi(q)|\J^Q_\C(q)|}{q}
	+O_{\C,\varepsilon}\big(Q^{-1/2+\varepsilon}\big)\,.
  \end{equation}
To calculate the sum in \eqref{eqpaispe}, we employ the Euler-MacLaurin
formula, noticing the fact that $|\J^Q_\C(q)|$, as a function of $q$, is
piecewise continuous differentiable on $[0,1]$. We obtain
  \begin{equation}\label{eqcinspe}
	\sum_{q=1}^Q \frac{\varphi(q)|\J^Q_\C(q)|}{q}
	=\frac{1}{\zeta(2)}\int\limits_{1}^Q |\J^Q_\C(q)|\,dq
	+O\bigg(\Big(\sup_{1\le q\le Q}|\J^Q_\C(q)|+
		\int\limits_{1}^Q 
		\frac{\partial}{\partial q}|\J^Q_\C(q)|\,dq \Big)\log Q\bigg)\,.
  \end{equation}
The size of the error term is estimated observing, firstly, that
$|\J^Q_\C(q)|\le Q$. Secondly, by the definition of $\J^Q_\C(q)$ it
follows that there exists a partition of $[1,Q]$ in finitely many
intervals with the property that the cardinality of $\J^Q_\C(q)$ is
monotonic on each of them. Therefore 
  \begin{equation}\label{eqsaispe}
	\int\limits_{1}^Q 
		\frac{\partial}{\partial q}|\J^Q_\C(q)|\,dq 
	= O_\C\big(Q\big)\,.
  \end{equation}
Then, forgathering \eqref{eqcinspe}, \eqref{eqsaispe}, \eqref{eqmax} in
\eqref{eqpaispe} and using again the estimate
$N_\I(Q)=3|\I|Q^2/\pi^2+O(Q\log Q)$, we obtain 
  \begin{equation*}
    \begin{split} 	
    \mu_{2,h}^{Q,\I}(\C)
	&=
	\frac{6|\I|}{\pi^2N_\I(Q)} \int\limits_{1}^Q |\J^Q_\C(q)|\,dq
		+O_{\C,\varepsilon}\big(Q^{-1/2+\varepsilon}\big) \\
	&=2\Area\big(\OO(Q)\big)
		+O_{\C,\varepsilon}\big(Q^{-1/2+\varepsilon}\big)\,.
    \end{split}
  \end{equation*}
This concludes the proof of the theorem.

\section{The Support of the Limiting Measure}\label{secSupport}

For $h=1$, we have $\D_{2,1}=[6/\pi^2,\infty)$.
For $h\ge 2$, by Theorem~\ref{Theorem1},
it follows that $\D_{2,h}$ is a countable union  
of hyper-surfaces in $[6/\pi^2,\infty)^h$.

The support $\D_{2,h}$ has some striking features. 
Let us see them in the case $h=2$.
We write ${\bf k}=(k,l)$ and observe that
\begin{equation*}
    \begin{split}
     \T_{k,l}=\Big\{(x,y)\in\T_k\colon\
              \frac{1+(l+1)x}{k(l+1)-1}< y\le \frac{1+lx}{kl-1}\Big\}\,.
    \end{split}
  \end{equation*}
Roughly speaking, by definition we find that 
$\T_{k}$ corresponds to the set of 3-tuples
$(\ga',\ga'',\ga''')$ of consecutive elements of $\FQ$ with the property that
$\Delta(\ga',\ga''')=k$. Similarly, $\T_{k,l}$ corresponds to the set of 4-tuples
$(\ga',\ga'',\ga''',\ga^{iv})$ of consecutive elements of $\FQ$
with the property that $\Delta(\ga',\ga''')=k$ and $\Delta(\ga'',\ga^{iv})=l$.
We remark that $\T_{1,1}=\emptyset$, and also $\T_{k,l}=\emptyset$ whenever both $k$
and $l$ are $\ge 2$, except when
$(k,l)\in \big\{(2,2); (2,3); (2,4); (3,2); (4,2)\big\}$. 
Notice that the symmetry of the Farey sequence of order $Q$ with respect to $1/2$ 
produces a sort of balance between the polygons $\T_{k,l}$ and $\T_{l,k}$.

Then the support $\D_{2,h}$ is the closure of the image of the
function $\Phi_{2,2}$, which can be written as
  \begin{equation*}
    \begin{split}
     \Phi_{2,2}(x,y)=\frac 3{\pi^2}\Big(\frac k{xz},\,\frac l{yt}\Big)\, ,
    \end{split}
  \end{equation*}
in which  $z=x-ky$, $t=y-lt$, for $(x,y)\in\T_{k,l}$ .
A tedious, but elementary, computation allows us to find precisely the boundaries of
$\Phi_{2,2}(\T_{k,l})$. The image obtained is shown in Figure~\ref{Figure1} and the
equations are listed in Table~\ref{Table1}.
All the functions that produce the equations of the boundaries
of $\Phi_{2,2}(\T_{k,l})$ are either of the form
$\frac{3}{\pi^2}\cdot\frac{et}{a+bt+c\sqrt{t(t-d)}}$, 
with $t$ in a certain interval that might be
unbounded, or the symmetric with respect to $x=y$ of such a curve.
Here $a,b,c,d,e$ are integers.

We conclude by making a few remarks.
Firstly, we mention that $\Phi_{2,2}$ has a symmetrization effect, namely, it 
makes $\Phi_{2,2}(\T_{n,m})$ and to $\Phi_{2,2}(\T_{m,n})$ to be symmetric  
with respect to the first diagonal $y=x$, for any $m,n\ge 1$.
The diamond\footnote{Remark that the edges of $\Phi_{2,2}(\T_{2,2})$ 
are close to being, but are not exactly straight lines. The same
applies for the edges of the diamonds in the tail.} $\Phi_{2,2}(\T_{2,2})$  
is the single nonempty domain
$\Phi_{2,2}(\T_{k,l})$ that has $y=x$ as axis of symmetry.
The top of the beak of the swallow $\D_{2,h}$ has coordinates
$(6/\pi^2,6/\pi^2)$. The asymptotes of the wings are $y=6/\pi^2$ and
$x=6/\pi^2$. The highest density is on a region situated in the neck,
where many components of the swallow overlap partially or completely.

\medskip
Table~\ref{Table1} below lists all the equations of the boundaries of 
$\Phi_{2,2}(\T_{k,l})$.
In the head of the table $MN$ represents an edge of $\T_{k,l}$ (listed in
counterclockwise order, starting either from the East or from the North side) and
$g_{_{MN}}(t)$ is a parametrization of $\Phi_{2,2}(MN)$.

\vspace*{0.1cm}
\bigskip
\small
\setlongtables
\begin{longtable}{|@{}D@{}|@{}C@{}|@{}C@{}|@{}C@{}|}
\caption{The edges of $\D_{2,2}$.}\label{Table1}\\ \hline
k,l & MN & \frac {\pi^2}{3}g_{_{MN}}(t) & \text{the\  domain\ of\ \ } t
         \\ 
\hhline{|====|}
\endfirsthead
\multicolumn{3}{l}{\small\sl continued from previous page}\\ 
\hline
k,l & MN & \frac {\pi^2}{3}g_{_{MN}}(t) & \text{the\  domain\ of\ \ } t 
        \\ 
\hhline{|====|}
\endhead
\multicolumn{4}{r}{\small\sl continued on next page} \\ 
\endfoot
\endlastfoot
1,2 & (\frac 13,1);\ (0,1) &\frac{2t}{\sqrt{t(t-4)}_{_{_{}}}} & \frac 92 \le t\le \infty \\
\hline
1,2 & (0,1);\ (\frac 15,\frac 45) &\frac{16t}{-12+3t+5\sqrt{t(t-8)}_{_{_{}}}} & \frac {25}3 \le t\le \infty \\
\hline
1,2 & (\frac 15,\frac 45);\ (\frac 13,1) &\frac{16t}{-12-3t+5\sqrt{t(t+8)}_{_{_{}}}} &
\frac 92 \le t\le \frac{25}{3} \\ 
\hhline{|====|}
1,3 & (\frac 12,1);\ (\frac 13,1) &\frac{6t}{t+3\sqrt{t(t-4)}\s} & 4 \le t\le \frac 92 \\
\hline
1,3 & (\frac 13,1);\ (\frac 15,\frac 45) &\frac{12t}{-t+3\sqrt{t(t+8)}\s} &
\frac 92 \le t\le \frac {25}{3} \\
\hline
1,3 & (\frac 15,\frac 45);\ (\frac 14,\frac 34)
&\frac{24t}{-20+7t+9\sqrt{t(t-8)}\s} &
8 \le t\le \frac {25}{3} \\
\hline
1,3 & (\frac 14,\frac 34);\ (\frac 27,\frac 57)
&\frac{24t}{-20+7t-9\sqrt{t(t-8)}\s} &
8 \le t\le \frac {49}{6} \\
\hline 
1,3 & (\frac 27,\frac 57);\ (\frac 12,1) &\frac{54t}{-24-7t+11\sqrt{t(t+12)}\s} &
4 \le t\le \frac {49}{6} \\
\hhline{|====|}
1,4 & (\frac 35,1);\ (\frac 12,1) &\frac{4t}{t-2\sqrt{t(t-4)}\s} & 
\frac{25}{6} \le t\le  4 \\
\hline
1,4 & (\frac 12,1);\ (\frac 27,\frac 57) &\frac{12t}{-t+2\sqrt{t(t+12)}\s} &
4 \le t\le \frac {49}{6} \\
\hline
1,4 & (\frac 27,\frac 57);\ (\frac 13,\frac 23)
&\frac{32t}{-28+11t-13\sqrt{t(t-8)}\s} &
\frac{49}{6} \le t\le 9 \\
\hline
1,4 & (\frac 13,\frac 23);\ (\frac 35,1)
&\frac{128t}{-40-13t+19\sqrt{t(t+16)}\s} &
\frac{25}{6} \le t\le 9 \\
\hhline{|====|}
1,\,  l\ge 5 & (\frac{l-1}{l+1},1);\ (\frac{l-2}{l},1) &
\frac{2lt}{(l-2)t-l\sqrt{t(t-4)}\s} & 
\frac{l^2}{2(l-2)} \le t\le  \frac{(l+1)^2}{2(l-1)} \\
\hline
1,\, l\ge 5 & (\frac{l-2}{l},1);\ (\frac{l-3}{l+1},\frac{l-1}{l+1}) &
\frac{2l(l-1)t}{(2-l)t+l\sqrt{t(t+4l-4)}\s} &
\frac{l^2}{2(l-2)} \le t\le  \frac{(l+1)^2}{2(l-3)} \\
\hline
1,\,  l\ge 7 & (\frac{l-3}{l+1},\frac{l-1}{l+1});\ (\frac{l-2}{l+2},\frac{l}{l+2}) &
\frac{8lt}{4+4l+(l-5)t-(l+3)\sqrt{t(t-8)}\s} &
\frac{(l+1)^2}{2(l-3)} \le t\le  \frac{(l+2)^2}{2(l-2)} \\
\hline
\,1,\,  l=5,6\ \, & (\frac{l-3}{l+1},\frac{l-1}{l+1});\ (\frac{l-2}{l+2},\frac{l}{l+2}) &
\frac{8lt}{4+4l+(l-5)t+(l+3)\sqrt{t(t-8)}\s} &
\frac{(l+1)^2}{2(l-3)} \le t\le  \frac{(l+2)^2}{2(l-2)} \\
\hline
1,\, l\ge 5 & (\frac{l-2}{l+2},\frac{l}{l+2});\ (\frac{l-1}{l+1},1) &
\frac{4l^3t^2}{\big((1-2l)t+\sqrt{t(t+4l)}\big)\big((l-1)t-(l+1)\sqrt{t(t+4l)}\big)\s}&
\frac{(l+1)^2}{2(l-1)} \le t\le  \frac{(l+2)^2}{2(l-2)} \\
\hhline{|====|}
2,1 & (1,1);\ (\frac 13,\frac 23) &\frac{4t^2}{(t+2)(t-2)\s} &
2 \le t\le  6 \\
\hline
2,1 & (\frac 13,\frac 23);\ (\frac 25,\frac 35) &
\frac{9t}{-12+4t-5\sqrt{t(t-6)}\s} &
6 \le t\le \frac{25}{4} \\
\hline
2,1 & (\frac 25,\frac 35);\ (1,1) &\frac{9t}{-12-4t+5\sqrt{t(t+6)}\s} &
2 \le t\le \frac{25}{4} \\
\hhline{|====|}
2,2 & (1,\frac 45);\ (1,1) &\frac{-8t^2}{(t+2)(t-6)\s} & 
2 \le t\le  \frac{10}{3} \\
\hline
2,2 & (1,1);\ (\frac 25,\frac 35) &\frac{6t}{-t+2\sqrt{t(t+6)}\s} &
2 \le t\le \frac {25}{4} \\
\hline
2,2 & (\frac 25,\frac 35);\ (\frac 12,\frac 12) &
\frac{18t}{-30+13t-14\sqrt{t(t-6)}\s} &
\frac{25}{4} \le t\le 8 \\
\hline
2,2 & (\frac 12,\frac 12);\ (1,\frac 45)
&\frac{50t}{-30-11t+14\sqrt{t(t+10)}\s} &
\frac{10}{3} \le t\le 8 \\
\hhline{|====|}
2,3 & (1,\frac 57);\ (1,\frac 45) &\frac{-12t^2}{(t+2)(t-10)\s} &
\frac{10}{3} \le t\le  \frac{14}{3} \\
\hline
2,3 & (1,\frac 45);\ (\frac 12,\frac 12) &\frac{30t}{-4t+6\sqrt{t(t+10)}\s} &
\frac{10}{3} \le t\le 8 \\
\hline
2,3 & (\frac 12,\frac 12);\ (\frac 45,\frac 35) &
\frac{27t}{24-2t+7\sqrt{t(t-6)}\s} &
8 \le t\le \frac{25}{4} \\
\hline
2,3 & (\frac 45,\frac 35);\ (1,\frac 57)
&\frac{147t}{-56-22t+27\sqrt{t(t+14)}\s} &
\frac{14}{3} \le t\le \frac{25}{4} \\
\hhline{|====|}
2,4 & (1,\frac 23);\ (1,\frac 57) &\frac{-16t^2}{(t+2)(t-14)\s} &
\frac{14}{3} \le t\le  6 \\
\hline
2,4 & (1,\frac 57);\ (\frac 45,\frac 35) &
\frac{28t}{-3t+4\sqrt{t(t+14)}\s} &
\frac{14}{3} \le t\le \frac{25}{4} \\
\hline
2,4 & (\frac 45,\frac 35);\ (1,\frac 23) & \frac{36t}{30-t+8\sqrt{t(t-6)}\s} &
6 \le t\le \frac{25}{4} \\
\hhline{|====|}
3,1 & (1,\frac 35);\ (1,\frac 23) & \frac{-9t^2}{(t+3)(t-6)\s} & 
3 \le t\le  \frac{15}{4} \\
\hline
3,1 & (1,\frac 23);\ (\frac 12,\frac 12) & \frac{9t^2}{(t+3)(2t-3)\s} &
3 \le t\le 6 \\
\hline
3,1 & (\frac 12,\frac 12);\ (\frac 47,\frac 37) &
\frac{32t}{-72+31t-11\sqrt{3t(3t-16)}\s} &
6 \le t\le \frac{147}{20} \\
\hline
3,1 & (\frac 47,\frac 37);\ (1,\frac 35) & 
\frac{50t}{-60-23t+9\sqrt{3t(3t+20)}\s} &
\frac{15}{4} \le t\le \frac{147}{20} \\
\hhline{|====|}
3,2 & (1,\frac 12);\ (1,\frac 35) & \frac{-18t^2}{(t+3)(t-15)\s} &
\frac{15}{4} \le t\le  6 \\
\hline
3,2 & (1,\frac 35);\ (\frac 47,\frac 37) & \frac{10t}{-2t+\sqrt{t(9t+60)}\s} &
\frac{15}{4} \le t\le \frac{147}{20} \\
\hline
3,2 & (\frac 47,\frac 37);\ (\frac 35,\frac 25) &
\frac{64t}{-168+79t-27\sqrt{3t(3t-16)}\s} &
\frac{147}{20} \le t\le \frac{25}{3} \\
\hline
3,2 & (\frac 35,\frac 25);\ (1,\frac 12) &
\frac{64t}{-72-11t+7\sqrt{3t(3t-16)}\s} &
6 \le t\le \frac{25}{3} \\
\hhline{|====|}
4,1 & (1,\frac 37);\ (1,\frac 12) & \frac{-16t^2}{(t+4)(t-12)\s} &
4 \le t\le  \frac{28}{5} \\
\hline
4,1 & (1,\frac 12);\ (\frac 35,\frac 25) & \frac{16t^2}{(t+4)(3t-4)\s} &
4 \le t\le \frac{20}{3} \\
\hline
4,1 & (\frac 35,\frac 25);\ (\frac 23,\frac 13) &
\frac{25t}{-80+37t-38\sqrt{t(t-5)}\s} &
9 \le t\le \frac{20}{3} \\
\hline
4,1 & (\frac 23,\frac 13);\ (1,\frac 37) &
\frac{49t}{-56-23t+26\sqrt{t(t+7)}\s} &
\frac{28}{5} \le t\le 9 \\
\hhline{|====|}
4,2 & (1,\frac 25);\ (1,\frac 37) & \frac{-32t^2}{(t+4)(t-28)\s} &
\frac{28}{5} \le t\le  \frac{20}{3} \\
\hline
4,2 & (1,\frac 37);\ (\frac 23,\frac 13) & \frac{14t}{-3t+4\sqrt{t(t+7)}\s} &
\frac{28}{5} \le t\le  9 \\
\hline
4,2 & (\frac 23,\frac 13);\ (1,\frac 25) & \frac{50t}{60-9t+16\sqrt{t(t-5)}\s} &
\frac{20}{3} \le t\le 9 \\
\hhline{|====|}
k\ge 5,\ 1 & (1,\frac {2}{k+1});\ (1,\frac 2k)
&\frac{-k^2t^2}{(t-k^2+k)(t+k)\s} & 
k \le t\le  \frac{k(k+1)}{k-1} \\
\hline
k\ge 5,\ 1 & (1,\frac 2k);\ (\frac{k-1}{k+1},\frac 2{k+1})
&\frac{k^2t^2}{(t+k)((k-1)t-k)\s} &
k \le t\le \frac{k(k+1)}{k-1} \\
\hline
k\ge 5,\ 1 & (\frac{k-1}{k+1},\frac{2}{k+1});\
(\frac{k}{k+2},\frac{2}{k+2}) &
\frac{4(k+1)^2t^2}{\big((k+2)t-\sqrt{kt(kt-4k-4)}\big)\big((k^2-2)t-k\sqrt{kt(kt-4k-4)}\big)\s}
&
\frac{k(k+1)}{k-1} \le t\le \frac{(k+2)^2}{k} \\
\hline
k\ge 5,\ 1 & (\frac{k}{k+2},\frac{2}{k+2});\ (1,\frac{2}{k+1})
& 
\frac{2(k+1)^2t^2}{\big(t-\sqrt{kt(kt-4k-4)}\big)\big((-k-2)t+\sqrt{kt(kt-4k-4)}\big)\s}
&
\frac{k(k+1)}{k-1} \le t\le \frac{(k+2)^2}{k} \\
\hline
\end{longtable}
\normalsize
\bibliographystyle{amsplain}

\end{document}